\title{Slicing convex sets and measures by a hyperplane}

\author{
{\sc Imre B\'ar\'any}
\\
  {\footnotesize R\'enyi Institute of Mathematics}\\[-1.5mm]
  {\footnotesize Hungarian Academy of Sciences}\\[-1.5mm]
  {\footnotesize POBox 127, 1364 Budapest, Hungary}\\[-1.5mm]
   {\footnotesize e-mail: {\tt barany@renyi.hu} and}\\[-1.5mm]
  {\footnotesize Department of Mathematics}\\[-1.5mm]
  {\footnotesize University College London}\\[-1.5mm]
  {\footnotesize Gower Street, London WC1E 6BT}\\[-1.5mm]
  {\footnotesize England},\\
{\sc Alfredo Hubard}\\
{\footnotesize Courant Institute of Mathematics of NYU}\\[-1.5mm]
{\footnotesize 251 Mercer Street 10012 NY}\\[-1.5mm]
{\footnotesize e-mail: {\tt hubard@cims.nyu.edu}}\\
{\sc Jes\'us Jer\'onimo}\\
{\footnotesize Centro de Investigaci\'on en Matem\'aticas A.C.}\\[-1.5mm]
{\footnotesize Apdo. Postal 402, Guanajuato, M\'exico}\\[-1.5mm]
{\footnotesize e-mail: {\tt jeronimo@cimat.mx}}}

\documentclass[12pt]{article}
\usepackage[pctex32]{graphics}
\usepackage{amsmath}
\usepackage{pst-all}
\usepackage{amscd}
\usepackage{amsfonts}

\newtheorem{teo}{Theorem}
\newtheorem{lem}{Lemma}
\newtheorem{prop}{Proposition}
\newtheorem{obs}{Observation}
\newtheorem{cor}{Corollary}
\newtheorem{defi}{Definition}
\newtheorem{prob}{Problema}
\newtheorem{sol}{Proof}
\newcommand{\bteo}{\begin{teo}\noindent {\bf }~}
\newcommand{\eteo}{\end{teo}}
\newcommand{\blem}{\begin{lem}\noindent {\bf }~}
\newcommand{\elem}{\end{lem}}
\newcommand{\bcor}{\begin{cor}\noindent {\bf }~}
\newcommand{\ecor}{\end{cor}}
\newcommand{\bprop}{\begin{prop}\noindent {\bf }~}
\newcommand{\eprop}{\end{prop}}
\newcommand{\bobs}{\begin{obs}\noindent {\bf }~}
\newcommand{\eobs}{\end{obs}}
\newcommand{\bprob}{\begin{prob}\noindent {\bf }~}
\newcommand{\eprob}{\end{prob}}
\newcommand{\bdefi}{\begin{defi}\noindent {\bf }~}
\newcommand{\edefi}{\end{defi}}
\newcommand{\bsol}{\begin{sol}\noindent {\bf }~}
\newcommand{\esol}{\end{sol}}

\newcommand{\al}{\alpha}
\newcommand{\be}{\beta}

\newcommand{\xx}{\overline x}
\newcommand{\la}{\lambda}
\newcommand{\ga}{\gamma}

\newcommand{\conv}{{\mathop {\rm conv} \nolimits}}

\newcommand{\Vol}{{\mathop {\rm Vol\,} \nolimits}}
\newcommand{\aff}{{\mathop {\rm aff} \nolimits}}

\newcommand{\R}{{\mathbb{R}^d}}
\newcommand{\FF}{{\mathcal F}}

\newcommand{\fin}{ $\Box$}

\begin{document}
\maketitle

\begin{abstract} Given convex bodies $K_1,\dots,K_d$ in $\R$ and numbers
$\al_1,\dots,\al_d \in [0,1]$, we give a sufficient condition for existence 
and
uniqueness of an (oriented) halfspace $H$ with $\Vol (H \cap K_i)=\al_i\cdot 
\Vol K_i$
for every $i$. The result is extended from convex bodies to measures.
\end{abstract}

\section{Transversal spheres}

A well known result in elementary geometry states that there is a unique 
sphere which
contains a given set of $d+1$ points in general position in $\R.$ A similar 
thing happens
with $d$-pointed sets and hyperplanes. What happens if we consider convex 
bodies instead
of points?

These questions are the main motivation for the present paper. The first 
result in this
direction is due to H. Kramer and A.B. N\'emeth \cite{kram}. They used the 
following,
very natural definition.

A family $\FF$ of connected sets in $\R$ is said to be {\sl well separated}, 
if for any
$k \le d+1$ distinct elements, $K_1, \ldots ,K_k$, of $\FF$ and for any 
choice of points
$x_i \in \conv K_i$, the set $\aff \,\{x_1,\dots,x_k\}$ is a $(k-1)$-dimensional 
flat. Here
$[k]$ stands for the set $\{1,2,\dots,k\}$. It is well known (cf. 
\cite{good} and
\cite{bar}), and also easy to check the following.

\begin{prop}\label{equiv} Assume $\FF=\{K_1,\dots,K_n\}$ is a family of 
connected sets in $\R$.
The following conditions are equivalent:
\begin{enumerate}
\item The family $\FF$ is well separated.
\item For every pair of disjoint sets $I,J \subset [n]$ with $|I|+|J| \le 
d+1$, there is
a hyperplane separating the sets $K_i,\; i \in I$ from the sets $K_j,\; j 
\in J$.
\end{enumerate}
\end{prop}

By an elegant application of Brouwer's fixed point theorem, Kramer and 
N\'emeth proved
the following:\newline

\textbf{Theorem KN.} \emph{Let $\FF$ be a well separated family of
$d+1$ compact convex sets in $\R$. Then there exists a unique
Euclidean ball which touches each set and whose interior is
disjoint from each member of $\FF$.}\newline

Denote by $B(x,r)$, resp. $S(x,r)$, the Euclidean ball and sphere
of radius $r$ and center $x$. We say that the sphere $S(x,r)$
\emph{supports} a compact set $K$ if $S(x,r)\cap K\neq\emptyset$
and either $K\subset B(x,r)$ or $K\cap \text{int}\,
B(x,r)=\emptyset.$ This definition is due to V. Klee, T. Lewis,
and B. Von Hohenbalken \cite{klee1}. They proved the
following:\newline

\textbf{Theorem KLH.} \emph{Let $\FF=\{K_1,K_2,\ldots ,K_{d+1}\}$ be a well 
separated
family of compact convex sets in $\R$, and let $I,J$ be a partition of 
$[d+1]$. Then
there is a unique Euclidean sphere $S(x,r)$ that supports each element of 
$\FF$ in such a
way that $K_i \subset B(x,r)$ for each $i\in I$ and $K_j \cap \text{int\;}
B(x,r)=\emptyset$ for each $j \in J$.}\newline

The case $I=\emptyset$ corresponds to Theorem KN. We are going generalize 
these results.
Let $Q^d=[0,1]^d$ denote the unit cube of $\R$. Given a well separated family 
$\FF$ of convex
sets in $\R$, a sphere $S(x,r)$ is said to be {\sl transversal} to $\FF$ if 
it intersects
every element of $\FF$. Finally, a {\sl convex body} in $\R$ is a convex 
compact set with
nonempty interior.

\begin{teo}\label{sphere}Let $\FF=\{K_1,\ldots ,K_{d+1}\}$ be a well 
separated family of
convex bodies in $\R$, and let $\al=(\alpha_1,\ldots ,\alpha_{d+1})\in 
Q^{d+1}$. Then
there exists a unique transversal Euclidean sphere $S(x,r)$ such that 
$\Vol(B(x,r)\cap
K_i)=\alpha_i\cdot \Vol(K_i)$ for every $i\in [d+1].$
\end{teo}

\textbf{Remark 1.} The transversality of $S(x,r)$ only matters when $\al_i$ 
is equal to
$0$ or $1$; otherwise the condition $\Vol(B(x,r)\cap K_i)=\alpha_i\cdot 
\Vol(K_i)$ plus
convexity guarantees that $S(x,r)$ intersects $K_i$.

\section{Transversal hyperplanes and halfspaces}

In a similar direction, S.E. Cappell, J.E. Goodman, J. Pach, R.
Pollack, M. Sharir, and R. Wenger \cite{capp} proved an analogous
theorem for the case of supporting hyperplanes, which can be seen
as spheres of infinite radius. Given a family $\FF$ of sets in
$\R$, a hyperplane will be called {\sl transversal} to $\FF$ if it
intersects each member of $\FF$. The following result is a special
case of Theorem 3 of Cappell et al. \cite{capp} (cf \cite{bisz} as
well):

\textbf{Theorem C.} \emph{Let $\FF=\{K_1,\ldots ,K_d\}$ be a well separated 
family of
compact convex sets in $\R$ with a partition $I,J$ of the index set $[d]$. 
Then there are
exactly two hyperplanes, $H_1$ and $H_2$, transversal to $\FF$ such that 
both $H_1$ and
$H_2$ have all $K_i$ ($i\in I$) on one side and all $K_j$ ($j\in J$) on the 
other side.}

Theorem C was also proved by Klee et al. \cite{klee2} using
Kakutani's extension of Brouwer's fixed point theorem. We are
going to formulate this theorem in a slightly different way, more
suitable for our purposes. So, we need to introduce new notation
and terminology.

A halfspace $H$ in $\R$ can be specified by its outer unit normal
vector, $v$, and by the signed distance, $t \in \mathbb R$, of its
bounding hyperplane from the origin. Thus, there is a one-to-one
correspondence between halfspaces of $\R$ and pairs $(v,t)\in
S^{d-1}\times \mathbb{R}$. We denote the halfspace $\{x\in \R:
\langle x,v\rangle \le t\}$ by $H(v\le t)$. Analogously we write
$H(v=t)=\{x\in \R: \langle x,v\rangle = t\}$, which is the
bounding hyperplane of $H(v\le t)$. Furthermore, given a set $K\subset
\R$, a unit vector $v$ and a scalar $t$, we denote the set
$H(v=t)\cap K$ by $K(v=t)$, analogously $K(v\leq t)=H(v\leq t)\cap
K$.

Suppose next that $\FF=\{K_1,\dots,K_d\}$ is a well separated family of 
convex sets in
$\R$. Assume $a_1 \in K_1,\dots,a_d \in K_d$. The unit normal vectors to the 
unique
transversal hyperplane containing these points are $v$ and $-v$. We want to 
make the
choice between $v$ and $-v$ unique and depend only on $\FF$. We first make 
it depend 
on $a_1,\dots,a_d$. Define $v=v(a_1,\dots,a_d)$ as the (unique) unit normal 
vector to
$\aff\{a_1,\dots,a_d\}$ satisfying
\[
\det \begin{vmatrix}
v & a_1 & a_2 & \cdots &a_d\\
0 & 1    &  1  & \cdots & 1
\end{vmatrix} > 0,
\]
in other words, the points $v+a_1,a_1,a_2,\dots,a_d$, in this order, are the 
vertices of
a positively oriented $d$-dimensional simplex. Clearly, with $-v$ in place 
of $v$ the
determinant would be negative. This gives rise to the map 
$v\;:K\longrightarrow S^{d-1}$
where $K=K_1\times \dots \times K_d$. This definition seems to depend on the 
choice of
the $a_i$, but in fact, it does not. Write $H(v=t)=\aff\{a_1,\dots,a_d\}$.

\begin{prop}\label{uniq} Under the previous assumption, let $b_i \in 
K_i(v=t)$ for each $i$.
Then $v(a_1,\dots,a_d)=v(b_1,\dots,b_d)$
\end{prop}

{\bf Proof.} This is simple. The homotopy $(1-\la)a_i+\la b_i$ ($\la \in 
[0,1]$) moves
the $a_i$ to the $b_i$ continuously, and keeps $(1-\la)a_i+\la b_i$ in 
$K_i(v=t)$. The
affine hull of the moving points remains unchanged, and does not degenerate 
because $\FF$
is well separated. So their outer unit normal remains $v$ throughout the 
homotopy. \fin

The previous proposition is also mentioned by Klee et al.
\cite{klee2}. With this definition, a transversal hyperplane to
$\FF$ determines $v$ and $t$ uniquely. We call $H(v= t)$ {\sl a
positive transversal hyperplane to} $\FF$, and similarly, $H(v\le
t)$ is {\sl a positive transversal halfspace to} $\FF$.

\begin{teo}\label{halfsp} Let $\FF=\{K_1,\ldots ,K_d\}$ be a family of well
separated convex bodies in $\R$, and let $\al=(\alpha_1,\ldots ,\alpha_d) 
\in Q^d$. Then
there is a unique positive transversal halfspace, $H$, such that $\Vol(K_i 
\cap
H)=\alpha_i\cdot \Vol(K_i)$ for every $i\in [d].$
\end{teo}

Theorem C follows since the partition $I,J$ gives rise to $\al, \be \in Q^d$ 
via
$\al_k=1$ if $k\in I$, otherwise $\al_k=0$, and $\be_k=1$ if $k \in J$, 
otherwise
$\be_k=0$. By Theorem \ref{halfsp}, there are unique positive transversal 
halfspaces
$H(\al)$ and $H(\be)$ with the stated properties. Their bounding hyperplanes 
satisfy the
statement of Theorem $C$ and they are obviously distinct. We mention, 
however, that
Theorem C will be used in the proof of the unicity part of Theorem 
\ref{halfsp}.

{\bf Remark 2.} When all $\alpha_i=1/2$, the existence of such a halfspace 
is guaranteed
by Borsuk's theorem, even without the condition of convexity or $\FF$ being 
well
separated. (Connectivity of the sets implies that the halving hyperplane is 
a transversal
to $\FF$.) The case of general $\alpha_i$, however, needs some extra 
condition as the
following two examples show. If all $K_i$ are equal, then each oriented 
hyperplane
section cuts off the same amount from each $K_i$, so 
$\alpha_1=\dots=\alpha_d$ must hold.
The second example consists of $d$ concentric balls with different radii, 
and if the
radius of the first ball is very large compared to those of the others and 
$\alpha_1$ is
too small, then a hyperplane cutting off $\alpha_1$ fraction of the first 
ball is
disjoint from all other balls. Thus no hyperplane transversal exists that 
cuts off an
$\alpha_1$ fraction of the first set.
\smallskip

\textbf{Remark 3.} Cappell et al. prove, in fact, a much more
general theorem \cite{capp}. Namely, assume that $\FF$ is well
separated and consists of $k$ strictly convex sets, $k \in
\{2,\dots,d\}$ and let $I,J$ be a partition of $[k]$. Then the set
of all supporting hyperplanes separating the $K_i$ $(i \in I)$ from
the $K_j$ $(j \in J)$ is homeomorphic to the $(d-k)$-dimensional
sphere.

\section{Extension to measures}

Borsuk's theorem holds not only for volumes but more generally for
measures. Similarly, our Theorem 2 can and will be extended to {\sl
nice measures} that we are to define soon. We need a small piece of
notation.

Let $\mu$ be a finite measure on the Borel subsets of $\R$ and let $v \in 
S^{d-1}$ be a
unit vector. Define
\begin{eqnarray*}
t_0&=&t_0(v)= \inf \{t \in \mathbb{R}:\; \mu(H(v\le t))>0\},\\
t_1&=&t_1(v)= \sup \{t \in \mathbb{R}:\; \mu(H(v\le t)) <
\mu(\R)\}.
\end{eqnarray*}
Note that $t_0=-\infty$ and $t_1=\infty$ are possible.

Let $H(s_0\le v\le s_1)$ denote the closed slab between the
hyperplanes\linebreak $H(v=s_0)$ and $H(v=s_1)$. Define the set
$K$ by
\[
K= \bigcap_{v \in S^{d-1}}H(t_0(v) \le v\le t_1(v)).
\]
$K$ is called the {\sl support} of $\mu$. Note that $K$ is convex 
(obviously) and
$\mu(\R \setminus K)=0$.

\begin{defi}\label{def} The measure $\mu$ is called {\sl nice} if the 
following conditions are satisfied:
\begin{enumerate}
\item[\rm(i)] $t_0(v)$ and $t_1(v)$ are finite for every $v \in S^{d-1}$,
\item[\rm(ii)] $\mu(H(v=t))=0$ for every $v \in S^{d-1}$ and $ t \in 
\mathbb{R}$,
\item[\rm(iii)] $\mu(H(s_0\le v \le s_1))>0$ for every $v \in S^{d-1}$ and 
for every
$s_0,s_1$ satisfying $t_0(v)\le s_0<s_1\le t_1 (v)$.
\end{enumerate}
\end{defi}

If $\mu$ is a nice measure, then its support is full-dimensional since, by 
(ii), it is
not contained in any hyperplane.

The function $t \mapsto \mu(K(v\le t))$ is zero on the interval 
$(-\infty,t_0]$, is equal
to $\mu(K)$ on $[t_1,\infty)$, strictly increases on $[t_0,t_1]$, and, in 
view of (iii),
is continuous. Assume $\alpha \in [0,1]$. Then there is a unique $t \in 
[t_0,t_1]$ with
\[
\mu(K(v\le t))=\alpha \cdot \mu(K).
\]
Denote this unique $t$ by $g(v)$; this way we defined a map $g: S^{d-1} 
\longrightarrow
\mathbb{R}$. The following simple lemma is important and probably well 
known.

\begin{lem}\label{contin}For fixed $\al \in [0,1]$ the function $g$ is 
continuous.
\end{lem}

{\bf Proof.} When $\alpha =1$, $g(v)$ is the support functional of $K$, 
which is not only
continuous but convex (when extended to all $v \in \R$). Similarly, $g$ is 
continuous
when $\alpha=0$.

Assume now that $0<\alpha <1$. Let $v_0\in S^{d-1}$ be an arbitrary point. 
In order to
prove the continuity of $g$ at $v_0$ we show first that $K(v=g(v))$ and 
$K(v_0=g(v_0))$
have a point in common whenever $v$ it is close enough to $v_0.$

Obviously, $K(v_0 = g(v_0))$ is a $(d-1)$-dimensional convex set lying in 
the hyperplane
$H(v_0=g(v_0))$. Then, for every small enough neighbourhood of $v_0,$ and 
for each $v$ in
such a neighbourhood, the supporting hyperplane of $K$ with unit normal $v$ 
(and $-v$) is
also a supporting hyperplane of $K(v_0\ge g(v_0))$ (and $K(v_0\le g(v_0))$).

Assume $s_v \le S_v$ and let $H(v= s_v)$ and $H(v= S_v)$ be the two 
supporting
hyperplanes (with normal $v$) to $K(v_0=g(v_0))$ which is a 
$(d-1)$-dimensional convex
set. Since $K(v_0=g(v_0))$ is a $(d-1)$-dimensional convex set, condition 
(iii) implies
that $s_v < S_v$. It follows that
\[
K(v\le s_v) \subset K(v_0 \le g(v_0)) \subset K(v\le S_v),
\]
and so
\[
\mu (K(v \le s_v)) \le \mu (K(v_0 \le g(v_0)) \le \mu (K(v \le S_v)).
\]
As $\mu (K(v_0 \le g(v_0))=\alpha \cdot \mu(K)$, we have $s_v \le g(v) \le 
S_v$.
Consequently, $K(v=g(v))$ and $K(v_0=g(v_0))$ have a point, say $z=z(v)$, in 
common. This
$z(v)$ is not uniquely determined but that does not matter.

It is easy to finish the proof now. Clearly $g(v)=\langle v, z(v)\rangle$ 
and
$g(v_0)=\langle v_0, z(v)\rangle$ for all $v$ in a small neighbourhood of 
$v_0$. Assume
the sequence $v_n$ tends to $v_0$. We claim that every subsequence, 
$v_{n'}$, of $v_n$
contains a subsequence $v_{n''}$ such that $\lim g(v_{n''})=g(v_0)$, which 
evidently
implies the continuity of $g$ at $v_0$.

For the proof of this claim observe first that, since $K(v_0=g(v_0))$ is 
compact,
$z(v_{n'})$ contains a convergent subsequence $z(v_{n''})$ tending to $z_0$, 
say. Taking
limits gives $z_o \in K(v_0=g(v_0))$. Then $g(v_{n''})=\langle v_{n''}, 
z(v_{n''})\rangle
\to \langle v_0, z_0\rangle =g(v_0)$. \fin\medskip

Theorem \ref{halfsp} is extended to measures in the following way.

\begin{teo}\label{meas} Suppose $\mu_i$ is a nice measure on $\R$ with 
support $K_i$ for all $i \in
[d]$. Assume the family $\FF=\{K_1,\ldots ,K_d\}$ is well separated and let
$\al=(\alpha_1,\ldots ,\alpha_d)\in Q^d$. Then there is a unique positive 
transversal
halfspace, $H$, such that $\mu_i(K_i\cap H)=\alpha_i\cdot \mu_i(K_i),$ for 
every $i\in
[d]$.
\end{teo}

\begin{cor}\label{cor}Assume $\mu_i$ are finite measures on $\R$ satisfying
conditions (i) and (ii) of Definition \ref{def}. Let $K_i$ be the support of 
$\mu_i$ for
all $i \in [d]$. Suppose the family $\FF=\{K_1,\ldots ,K_d\}$ is well 
separated and let
$\al=(\alpha_1,\ldots ,\alpha_d) \in Q^d$. Then there is a positive 
transversal
halfspace, $H$, such that $\mu_i(K_i\cap H)=\alpha_i\cdot \mu_i(K_i),$ for 
every $i\in
[d]$.
\end{cor}

The corollary easily follows from Theorem \ref{meas}; we omit the simple 
details.

Theorem \ref{halfsp} is a special case of Theorem \ref{meas}: when $\mu_i$ 
is the
Lebesgue measure (or volume) restricted to the convex body $K_i$ for all  $i 
\in [d]$ and
the family $\FF$ is well separated. Also, Theorem C is a special case of 
Theorem
\ref{meas}: when $\mu_i$ and $K_i$ are the same as above, and, for a given 
partition
$I,J$ of $[d]$, we set $\alpha_i=1$ for $i\in I,$ and $\alpha_j=0$ for $j\in 
J.$ Theorem
\ref{sphere} follows from Theorem \ref{meas} via ``lifting to the 
paraboloid''. This is
explained in the last section.

\section{Proof of Theorem \ref{meas}.}

In the proof we will use Brouwer's fixed point theorem. We will define a 
continuous
mapping from a topological ball to itself, such that a fixed point of this 
map yields a
halfspace with the desired properties. Set $K=K_1\times \ldots \times K_d$. 
Given a point
$x=(x_1,\ldots ,x_d)\in K$ we consider the hyperplane $\aff\, \{x_1,\ldots 
,x_d\}$. Since
the family $\FF$ is well separated, this hyperplane is well defined for each 
$x\in K$.
Let $H(v\leq t)$ be the (unique) positive transversal halfspace whose 
bounding hyperplane
is $\aff\, \{x_1,\ldots ,x_d\}$.

In Section 2 we defined the map $v: K\longrightarrow S^{d-1}$ which is the 
properly
chosen unit normal to $\aff\{x_1,\ldots ,x_d\}.$ Clearly, this function is 
continuous.

We prove existence first. We start with the case when $\alpha_i\in (0,1)$ 
for every $i\in
[d]$. We turn to the remaining case later by constructing a suitable 
sequence of
halfspaces.

Let $g_i: S^{d-1}\longrightarrow \mathbb{R}$ be the function such that for 
each $v\in
S^{d-1},$ $g_i(v)$ is the real number for which $\mu_i (K_i(v \leq 
g_i(v)))=\alpha \cdot
\mu_i (K_i)$ for each $i\in [d].$ Each $g_i$ is a continuous function by 
Lemma
\ref{contin}. Let $h:S^{d-1}\longrightarrow K$ be the function sending $v 
\mapsto
(s_1,\ldots,s_d)$ where $s_i$ is the Steiner point of the 
$(d-1)$-dimensional section,
$K_i(v=g_i(v))$ for each $i\in [d]$. As is well known, the family of 
sections $K_i(v=t)$
depend continuously (according to the Hausdorff metric) on the corresponding 
family of
hyperplanes, $\{H(v=t)\}$ whenever every section is $(d-1)$-dimensional, 
which is
obviously the case because $\alpha_i \in (0,1)$. It is also well known that 
the function
that assigns to a compact convex set its Steiner point is continuous. Hence, 
$h$ is a
continuous function.

It follows that $$f:=h \circ v: K \longrightarrow K$$ is a continuous 
function. As $K$ is
a compact convex set in $\R\times\ldots\times\R$ Brouwer's fixed point 
theorem implies
the existence of a point $x\in K$ such that $f(x) =x.$ Consider a fixed 
point,
$x=(x_1,\ldots ,x_d)$, of $f$. Then the halfspace $H(v\le t)$ whose bounding 
hyperplane
is $\aff\, \{x_1,\ldots ,x_d\}$ is a positive transversal halfspace to $\FF$ 
and it has
the required properties.

Next we prove existence for vectors $\alpha=(\alpha_1,\ldots ,\alpha_d)\in 
Q^d$ that may
have $0,1$ components as well. Consider the sequence $\{\alpha^n\}\subset 
Q^d$
$\alpha^n=(\alpha_1^n,\ldots,\alpha_d^n)$ (defined for every $n\geq 2$), 
such that for
every entry $\alpha_i=0$ we define $\alpha_i^n=\frac{1}{n},$ for every entry 
$\alpha_i=1$
we define $\alpha_i^n=1-\frac{1}{n},$ and for every entry $\alpha_i\notin 
\{0,1\}$ we
define $\alpha_i^n=\alpha_i.$ Also, for every $n\geq 2$ we consider the 
unique positive
transversal halfspace $H(v_n\le t_n)$ with $\mu_i (K_i(v_n \leq 
t_n))=\alpha_i^n \cdot
\mu_i (K_i),$ for each $i$.  The compactness of $K$ implies that the set of 
all possible
$(v,t)\in S^{d-1}\times \mathbb R$ such that the hyperplane $H(v=t)$ is 
transversal to
$\FF$ is compact. Thus there exists a convergent subsequence 
$\{(v_{n'},t_{n'})\}$ which
converges to a point $(v,t)\in S^{d-1}\times {\mathbb R}$. Clearly, $H(v\le 
t)$ is a
positive transversal halfspace to $\FF$ which satisfies $\mu_i(K_i(v\leq
t))=\alpha_i\cdot \mu_i(K_i)$ for every $i$.

Next comes uniqueness. We start with the $0,1$ case, that is, when
$\alpha=(\alpha_1,\ldots ,\alpha_d)$ with all $\al_i \in \{0,1\}$. Such an 
$\al$ defines
a $\be \in Q^d$ via $\be_i=1-\al_i$ for every $i$. By the previous existence 
proof there
is a unique positive transversal halfspace $H(v\le t)$ for $\al$ and another 
one $H(u\le
s)$ for $\be$. These two halfspaces are distinct, first because $u=v$ is 
impossible, and
second because of the following fact which implies that $u\ne -v$.

\begin{prop} For every pair of points $(a_1,\dots,a_d)$ and 
$(b_1,\dots,b_d)$ in $K$,
$v(a_1,\dots,a_d)$ and $-v(b_1,\dots,b_d)$ are distinct.\end{prop}

{\bf Proof.} Assume $v(a_1,\dots,a_d)= -v(b_1,\dots,b_d)$. Then the affine 
hulls of the
$a_i$ and the $b_i$ are parallel hyperplanes. We use the same homotopy as in 
the proof of
Proposition \ref{uniq}. As $\la$ moves from $0$ to $1$, the moving points 
$(1-\la)a_i+\la
b_i$ stay in $K_i$, and their affine hull remains parallel with 
$\aff\{a_1,\dots,a_d\}$.
So the outer normal remains unchanged throughout the homotopy. A 
contradiction. \fin

The condition $\mu_i( K_i(v \le t))= \al_i \mu_i(K_i)$ implies, in the given 
case, that
all $K_i$ ($i\in I$) are in $H(v\le t)$ and all $K_j$ $(j \in J)$ are in 
$H(v\ge t)$.
Thus $H(v=t)$ is a transversal hyperplane satisfying the conditions of 
Theorem C with
partition $I,J$ where $I=\{i \in [d]\;: \al_i=0\} \mbox{ and } J=\{j \in 
[d]\;:
\al_j=1\}$. The same way, $H(u=s)$ is a transversal hyperplane satisfying 
the conditions
of Theorem C with the same partition $J,I$.

The uniqueness of $H(v\le t)$ follows now easily. If we had two distinct 
positive
transversal halfspaces $H(v_1 \le t_1)$ and $H(v_2 \le t_2)$ for $\al$, then 
we would
have four distinct transversal hyperplanes with $K_i$ ($i \in I$) on one 
side and $K_j$
($j\in J$) on the other side, contradicting Theorem C.

Now we turn to uniqueness for general $\al$. Assume that there are two 
distinct positive
transversal halfspaces $H(v_1 \le t_1)$ and $H(v_2 \le t_2)$ for $\al$. 
Their bounding
hyperplanes cannot be parallel. Define $M=H(v_1 \le t_1) \cap H(v_2 \le 
t_2)$ and
$N=H(v_1 \ge t_1) \cap H(v_2 \ge t_2)$. The partition $I,J$ of the index set 
$[d]$ is
defined as follows: $i\in I$ if $M \cap \text{int}\,K_i\neq\emptyset$ and 
$j\in J$ if
$M\cap \text{int}\,K_j=\emptyset$. Set $K_i'=M \cap K_i$ for every $i\in I$ 
and $K_j'=N
\cap K_j$ for every $j\in J.$ Let $\FF '$ be the family consisting of all 
the convex
bodies $K_i'$ ($i \in I$) and $K_j'$ ($j \in J$). It is quite easy to see 
that no member
of $\FF '$ is empty. Moreover, $\FF'$ is evidently well separated. Given the 
partition
$I,J$, define $\ga$ by $\ga_i=1$ for $i \in I$ and $\ga_j=0$ for $j\in J$. 
Then there are
two transversal halfspaces (with respect to $\FF ')$, namely $H(v_k\le t_k)$ 
$k=1,2$
satisfying $\mu_i(K_i(v_k\le t_k))=\ga_i\mu_i(K_i)$ for every $i$. But every 
$\ga_i \in
\{0,1\}$ and we just established uniqueness in the $0,1$ case.

\fin\medskip

\section{Proof of Theorem 2}

We will use the well-known technique of lifting the problem from $\R$ to a 
paraboloid in
$\mathbb{R}^{d+1}$, and then apply Theorem 3.

In this section we change notation a little. A point in $\R$ is denoted by
$x=(x_1,\dots,x_d)$, a point in $\mathbb{R}^{d+1}$ is denoted by
$\xx=(x_1,\dots,x_d,x_{d+1})$. The projection of $\xx$ is 
$\pi(\xx)=(x_1,\dots,x_d)$, and
the lifting of $x$ is $\ell(x)=(x_1,\dots,x_d,|x|^2)$ where 
$|x|^2=x_1^2+\dots+x_d^2$.
Clearly, $\ell(x)$ is contained in the paraboloid
\[
P=\{\xx \in \mathbb{R}^{d+1}: \xx=(x_1,x_2,\ldots,x_d,|x|^2)\}.
\]
A set $K \subset \R$ lifts to $\ell(K)=\{\ell(x) \in P:\; x \in K\}$. Also,
$\pi(\ell(K))=K$.

A hyperplane is called {\sl non-vertical} if $\pi(H)=\R$. The lifting gives 
a bijective
relation between non-vertical hyperplanes in $\mathbb{R}^{d+1}$ 
(intersecting $P$) and
$(d-1)$-dimensional spheres in $\R$ in the following way. Assume $S=S(u,r)$ 
is the sphere
centered at $u$, with radius $r$ in $\R$. Of course, $\ell(S) \subset P$, 
but more
importantly,
\[
\ell(S)=P\cap H
\]
where $H$ is the hyperplane with equation $x_{d+1}=2\langle u,x 
\rangle+r^2-|u|^2$.
Conversely, given a non-vertical hyperplane $H$ with equation 
$x_{d+1}=2\langle u,x
\rangle +s$ where $s=r^2-|u|^2$ with some $r>0$,
\[
\pi(H \cap P)=S(u,r).
\]
As a first application of this lifting, here is a simple proof of a slightly 
stronger
version of Theorem KLH (we can replace the convexity assumption by 
connectedness).
Consider a family of $d+1$ well separated connected compact sets in $\R$ and 
a partition
of the sets into two classes. Lift the family into the paraboloid, and for 
each lifted
set, consider its convex hull. This gives a $(d+1)$-element family of convex 
bodies in
$\mathbb{R}^{d+1}$. The lifted family is well separated. This can be seen 
using
Proposition~\ref{equiv}: the lifting of the separating $(d-1)$-dimensional 
planes of the
original family yield (vertical) separating hyperplanes of the corresponding 
lifting.
Thus Theorem 3 applies to the lifted family (with the obviously induced 
partition) and
gives a hyperplane $H$ such $H \cap P$ projects onto a sphere $S$ in $\R$ 
satisfying the
requirements of Theorem~\ref{sphere}. We omit the straightforward detail.

We apply Theorem~\ref{meas} to the paraboloid lifting to obtain
Theorem~\ref{sphere}, in the same way. The family
$\FF=\{K_1,\dots,K_{d+1}\}$ lifts to the family
$\ell(\FF)=\{\ell(K_1),\dots,\ell(K_{d+1})\}$, and we define the
measures $\mu_i$ via
\[
\mu_i(C)=\Vol \pi(C \cap \ell(K_i)),
\]
where $C$ is a Borel subset of $\mathbb{R}^{d+1}$. Clearly, $\mu_i$ is 
finite and
$\ell(\FF)$ is well separated. Its support is $\conv \, \ell(K_i)$. It is 
easy to see
that $\mu_i$ is a nice measure by checking that it satisfies all three 
conditions.

Thus Theorem~\ref{meas} applies and guarantees the existence of a unique 
positive
transversal (to $\ell(\FF)$) halfspace $H \subset \mathbb{R}^{d+1}$ with 
$\mu_i(H\cap
\ell(K_i))=\alpha_i \cdot \mu_i(K_i)$ for each $i$. This translates to the 
ball $B=\pi(H
\cap P)$ and sphere $S=\pi(H^0 \cap P)$ (where $H^0$ is the bounding 
hyperplane of $H$)
as follows: $S$ is a transversal sphere of the family $\FF$ and $\Vol (B 
\cap
K_i)=\alpha_i \cdot \Vol K_i$. Unicity of $S$ follows readily. \fin

\section{Acknowledgments} The first author was partially supported by 
Hungarian National
Foundation Grants T 60427 and NK 62321. The second and third authors are 
grateful, for
support and hospitality, to the Department of Mathematics at University 
College London
where this paper was written. We also thank two anonymous referees for 
careful reading
and useful remarks and corrections.

\end{document}